\documentclass[oneside,12pt]{amsart}
\usepackage{amssymb, amscd}
\newcommand{\C}{{\mathbb C}}
\newcommand{\Q}{{\mathbb Q}}
\newcommand{\Z}{{\mathbb Z}}

\DeclareMathOperator{\im}{{im}}
\DeclareMathOperator{\id}{{Id}}
\DeclareMathOperator{\fix}{{Fix}}
\DeclareMathOperator{\ind}{{Ind}}
\DeclareMathOperator{\tor}{{tor}}

\theoremstyle{plain}
\newtheorem{theorem}{Theorem}[section]
\newtheorem{lemma}[theorem]{Lemma}

\newtheorem*{theorem*}{}
\newtheorem{proposition}[theorem]{Proposition}
\newtheorem{corollary}[theorem]{Corollary}

\theoremstyle{definition}

\newtheorem{definition}[theorem]{Definition}

\newtheorem{example}[theorem]{Example}
\newtheorem{remark}[theorem]{Remark}

\begin{document}\title{Fixity and
Free Group Actions on Products of Spheres}
\author{Alejandro
Adem$^*$}
\date{\today}
\address{Department of Mathematics\\        
University of Wisconsin\\         Madison WI
53706}
\email{adem@math.wisc.edu}
\author{James F.
Davis$^*$}
\address{Department of Mathematics\\	 Indiana University\\	
Bloomington IN 47405}
\email{jfdavis@indiana.edu}
\author{\"Ozg\"un
\"Unl\"u}\address{Department of Mathematics\\	 University of Wisconsin\\	
Madison WI 53706}
\email{unlu@math.wisc.edu}
\thanks {$^*$The first two authors were partially supported by
the NSF}

\begin{abstract}
A representation
$G\subset U(n)$ of degree $n$ has
\emph{fixity} equal to the smallest integer $f$ such that the induced action
of $G$ on
$U(n)/U(n-f-1)$ is free. Using bundle theory we show that if $G$ 
admits a representation of fixity one, then it acts freely and smoothly on
$\mathbb S^{2n-1}\times\mathbb S^{4n-5}$. 
We use this to prove that a finite $p$--group (for $p>3$)
acts freely and
smoothly on a product of two spheres if and only if it does not
contain $(\mathbb Z/p)^3$ as a subgroup.

We use propagation methods from surgery theory
to show that a representation of fixity $f<n-1$ gives rise to a free
action of $G$
on a product of $f+1$
spheres provided the order of $G$
is relatively prime to
$(n-1)!~$.  
We  give an
infinite collection of new examples of finite $p$--groups of rank
$r$ which act freely on a product of $r$ spheres, hence verifying a strong
form of a well-known conjecture for these groups. In addition we show that
groups of fixity two act freely on a finite complex
with the homotopy type of a product of three spheres. A number of
examples
are explicitly described.
\end{abstract}
\maketitle
\tableofcontents
\section{Introduction}
A well known result in topology is the characterization  of those finite groups
that can act freely on a sphere, namely groups with either cyclic or generalized
quaternion Sylow subgroups and such that every involution is central (see \cite{MTW}). 
Free linear actions on spheres are the most basic examples.
These can be constructed as follows: a subgroup $G$ of $U(n)$ 
acts on the homogeneous space $\mathbb S^{2n-1}\cong U(n)/U(n-1)$; if 
no conjugate of $G$ intersects $U(n-1)$ non--trivially,
then this gives rise to a free linear action on the sphere.

In this paper we consider the situation for products of spheres. 
Here the problem 
is much more complex, and in particular we still do not have a characterization of
those finite groups that can act freely on a product of {\sl two} spheres. 
A key new ingredient in our approach is the use
of the $G$ actions on the homogeneous spaces $U(n)/U(k)$; our view is that these
are also important building blocks for actions of groups of larger rank. 
On the algebraic side,
group theorists have studied representations such that the action of $G$ on
$U(n)/U(n-f-1)$ is free; the smallest such $f$ is called the 
\emph{fixity} of the representation. In particular there is an explicit
classification
of those finite $p$--groups having a representation of fixity $f<p$. 

Our main results are in the direction of propagating the natural free actions
on $U(n)/U(k)$ to free actions on an actual product of spheres. For low fixity
we can use explicit arguments involving equivariant vector bundles;
in the case of
fixity equal to one we have

\begin{theorem}
A subgroup $G\subset U(n)$ of fixity one acts freely and smoothly on
$X=\mathbb S^{2n-1}\times\mathbb S^{4n-5}$.
In particular if $G$ is any finite subgroup of $SU(3)$, then it will act freely 
and smoothly on $\mathbb S^5\times\mathbb S^7$.
\end{theorem}
\noindent From this we conclude that $A_5$, $SL_3(\mathbb F_2)$
and the triple cover $3A_6$ all act freely and smoothly on $\mathbb S^5\times \mathbb S^7$.

Using the explicit classification of rank two $p$--groups
(see \cite{Blackburn}) as well as the description of $p$--groups
of fixity equal to one in \cite{NOS} we see that our result can be
used to construct free actions for the exceptional (i.e. non--metacyclic)
$p$--groups on the
list. Combining this with the
well--known fact that metacyclic groups act freely on a product
of two linear spheres, we obtain

\begin{theorem}
Let $p>3$ be a prime. Then a $p$--group $P$ acts freely and
smoothly on $\mathbb S^m\times\mathbb S^n$ for some $m,n >0$
if and only if $P$ does not contain $\mathbb Z/p\times\mathbb Z/p
\times\mathbb Z/p$ as a subgroup.
\end{theorem}
\noindent This geometric result improves on the
homotopy theoretic version recently established in \cite{AS}.
Indeed this paper arose out of efforts to promote the results
there to actions on manifolds. Moreover we obtain all of
these actions explicitly, thus providing rank 2 models for 
group actions which may play a special role analogous to 
that of linear spheres. In contrast, the situation at the primes
$p=2, 3$ remains unresolved, reflecting the
complications in the corresponding group theory at these primes.
Recently \"Unl\"u \cite{Unlu} has shown that among the 396 groups of
order dividing $256$ and which have rank equal to two, there is exactly
\textbf{one} group which is not known to act freely and smoothly on a 
product of two spheres!

Using the methods developed in \cite{AS} we also obtain interesting results
for the case of fixity equal to two, namely

\begin{theorem}
If $G\subset U(n)$ is of fixity equal to two, then $G$ acts freely on a
finite complex $X\simeq\mathbb S^{2n-1}\times\mathbb S^{4n-5}
\times\mathbb S^M$ for some $M>0$.
In particular if $G\subset SU(4)$, then $G$ acts freely on a finite complex
$X\simeq \mathbb S^7\times\mathbb S^{11}\times\mathbb S^M$ for some $M>0$.
\end{theorem}
\noindent For example, this can be used to show that $Sp_4(\mathbb F_3)$ acts
freely on a finite complex $X\simeq \mathbb S^7\times\mathbb S^{11}
\times\mathbb S^M$ for some $M>0$.

For arbitrary fixity we must make use of methods from propagation theory,
involving homotopy theory and surgery. Observe that the Stiefel manifolds
$U(n)/U(n-f-1)$ have the cohomology of a product of $f+1$ spheres (see
\ref{cohomology}); our goal is to propagate this natural model to an action
on a product of spheres.
The main result in this paper is
the following

\begin{theorem}
Let $G$ denote a finite subgroup of $U(n)$ which
acts freely on
$U(n)/U(k)$ with $k \geq 1$. 
If the order of $G$ is prime to $(n-1)!~$, then $G$ acts freely, smoothly 
and homologically trivially on 
$\mathbb S^{2n-1}\times\mathbb S^{2n-3}\times\dots\times\mathbb S^{2k+1}$.
\end{theorem}
\noindent Note that if $G\subset SU(n)$, then $G$ will act freely
on $U(n)/U(1)$ (see \ref{SU(n)}), hence the theorem always applies.
More generally, if $G\subset U(n)$ and $(|G|,(n-1)!)=1$, our methods yield a free 
$G$--action
on a finite complex $X\simeq \mathbb S^{2n-1}\times\dots\times\mathbb S^3
\times\mathbb S^M$ for some $M>1$; we conjecture that an analogue
of our main theorem
should also hold in this case, but there are surgery--theoretic
difficulties to overcome which we hope to address in a subsequent
paper.

Applying the available characterization of low fixity $p$--groups, we
obtain interesting examples of group actions:
\begin{theorem}
Let $P$ denote a finite non--abelian $p$--group with cyclic center and
having an abelian maximal subgroup. If the rank of $P$ is $r<p$, then 
there exists a free and homologically trivial action of $P$ on
$M=\mathbb S^{2p-1}\times\mathbb S^{2p-3}\times\dots\times \mathbb S^{2(p-r)+1}$,
a product of $r$ spheres.
\end{theorem}

We should point out that this produces an infinite number of new examples
of free actions
by rank $r$ groups 
on a product of $r$ spheres. More generally it is conjectured that a 
rank $r$ finite group
will act freely on a finite complex $X$ having the homotopy type of a
product of $r$ spheres. A related conjecture is that every finite group
acts freely and homologically trivially on a product of spheres.
Although the condition
$(|G|, (n-1)!)=1$ is somewhat restrictive, our approach does yield a 
new method for approaching the conjectures mentioned above; most importantly
we have constructed many interesting \emph{geometric} actions.

We are grateful to R.L. Griess for providing 
information concerning the classification of finite
linear groups.

\section{Basic Definitions and Properties} In this
section
we will recall the notion of fixity for a complex representation of a finite
group $G$ (see \cite{NOS}) and relate it to properties of associated actions
on complex Stiefel manifolds $U(n)/U(k)$. 

First we introduce some notation.
Let $V$ denote a finite dimensional $\mathbb C G$--module.  For a subgroup
$H \subset G$ , we denote by $V^H$ the subspace of vectors
in $V$ fixed by all $h\in H$.
We  denote by 
$\langle g\rangle$ the subgroup generated by an element $g$ of $G$.

\begin{definition}
The
\emph{fixity} of a finite dimensional $\mathbb C G$--module $V$ is 
$$\fix_G(V) = \max_{g\in G} \{ \rm{dim}_{\mathbb C}~V^{\langle
g\rangle}\}$$
\end{definition}

Given a faithful complex representation $V$ of $G$, we can always
obtain an equivalent unitary representation. Hence in what follows
we will
restrict our attention to faithful unitary representations and the
associated embeddings $G\hookrightarrow U(n)$, where $U(n)$ denotes
the group of $n\times n$ unitary matrices and $n=\dim V$. Given such
an embedding and a closed subgroup $H\subset U(n)$, we have a natural 
$G$--action on the
homogeneous space $U(n)/H$.
In particular we will be interested in the subgroups $U(k)\subset U(n)$,
defined by
$$A\mapsto \left(\begin{array}{cc}A&0\\
0&I_{n-k}\\\end{array}\right)$$ 
where $I_{n-k}$ is the
$(n-k)\times (n-k)$ identity matrix.

Note that $\fix_G(V)$ can be expressed as
the maximum (for $g\in G$)
of the dimensions of the eigenspaces $\ker (g-I)$.

\begin{lemma}
$\fix_G(V) \le f$ if and only if the induced action of $G$ on $U(n)/U(n-f-1)$
is a free action.
\end{lemma}

\begin{proof}
If $\fix_G(V)\le f$, then
given any $1\ne g\in G$, $\ker (g-I)$ can be at most $f$--dimensional, hence $g$
cannot be conjugated into the subgroup $U(n-f-1)$ and so the action of $G$
on the homogeneous space $U(n)/U(n-f-1)$ must be free. Conversely if the
action is free, no element can be conjugated into $U(n-f-1)$, hence
$\fix_G(V)\le f$.
\end{proof}

\begin{corollary} 
A faithful $\mathbb C G$--
module $V$ has fixity $f$ if and only if $f$ is the smallest integer such that
the induced $G$--action on $U(n)/U(n-f-1)$ is
free.
\end{corollary}

We can now reformulate the notion of fixity.

\begin{definition} 
A faithful unitary representation $G \subset U(n)$
has {\em fixity $f$} if $f$ is the smallest integer so that the induced action
of $G$ on
$U(n)/U(n-f-1)$ is free. 
\end{definition}

 From this we derive an invariant associated to $G$.

\begin{definition}
For a finite group $G$, we define its \emph{fixity}, $\fix(G)$, as
the minimum value of $\fix_G(V)$, as $V$ ranges over all faithful,
finite dimensional complex representations of $G$.  
\end{definition}

The case of fixity zero coincides with that of linear space
forms. Indeed if an $n$--dimensional representation $V$ has fixity
zero then the induced action on $U(n)/U(n-1)$ is free; this can be
identified with the action on the sphere $S(V)$ of unit vectors in
$V$. Examples of higher fixity are less well--known from the point
of view of transformation groups, and they will provide building
blocks for new examples of group actions.

\begin{remark}\label{SU(n)}
Note that $SU(n)$ is a normal
subgroup of $U(n)$ which intersects $U(1)$ trivially.
Hence if $G\subset SU(n)$, it will act freely on $U(n)/U(1)$
and so the 
fixity is at most $n-2$.
\end{remark}

We now relate fixity to
another group--theoretic invariant.

\begin{definition}
For a finite group $G$ and a prime $p$, we define its \emph{$p$--rank} as 
$r_p(G)=\max~\{r~|~(\mathbb Z/p)^r\subset G\}$ and its \emph{rank} as
$r(G)=\max~\{r_p(G)~|~p~\rm{divides}~|G|\}$.
\end{definition}

\noindent The following is a basic result relating rank to fixity
(\cite{Sh}, Lemma 3.1):

\begin{proposition}
For any finite group $G$, $r(G)\le \fix(G)+1$. 
\end{proposition}

Note that this inequality may be strict. For example, if $G$ is a non-abelian split
extension of $\mathbb Z/p$ by $\mathbb Z/q$ where $p$ and $q$ are prime, then this is a rank one group which does not
have a fixed point free representation (i.e. $G$ is not a Frobenius complement, as it does
not satisfy the $pq$ condition); hence its fixity is greater than zero.  For us
the most interesting case occurs when $r(G)=\fix(G) +1$. 

Assume that
$G\subset U(n)$ has fixity equal to $f$. Then the $G$ action on $U(n)/U(n-f)$
is not free. In fact we have

\begin{proposition}
If $G\subset U(n)$ has fixity
$f$, then all of the isotropy subgroups for the $G$--action on $U(n)/U(n-f)$
have rank equal to one. 
\end{proposition}

\begin{proof}
Let $H$ denote an isotropy subgroup for the $G$--action on
$U(n)/U(n-f)$; this means that there exists a $g\in U(n)$ such that
$g^{-1}Hg\subset U(n-f)$. This subgroup then acts freely on the quotient
$U(n-f)/U(n-f-1)\cong \mathbb S^{2(n-f)-1}$ as otherwise a conjugate
of $G$ would intersect $U(n-f-1)$ non--trivially. Hence $H$ must be a group
of rank equal to one.
\end{proof}


\begin{remark}
More generally if we let
$M_k=U(n)/U(k)$, this defines a sequence of $G$--manifolds 
$\{M_0,M_1,\dots,M_n\}$ 
and equivariant maps $M_k\to
M_{k+1}$, for $k=0,\dots ,n$ where $M_0$ is a free $G$--space, $M_{n-1}=\mathbb S^{2n-
1}$ and $M_n=\{x_0\}$. As we go up this tower, the isotropy must
increase from rank zero to rank $r(G)$. Using an argument similar to the one above, 
we see that the rank can only increase by
one
at each stage. There are $n$ steps and a total increase by $r(G)$ must
happen. Hence we have a partition of $[0,n] \cap \Z$, given by integers
$0\le s_0<s_1<\dots < s_{r(G)}\le n$
such that the isotropy of $M_k$ has rank $t$ if $s_t\le k < s_{t+1}$.
The patterns which arise in this process seem to be an interesting
invariant of the representation; we shall make use of this in our
applications.
Note that when $n=r(G)$,
the rank must increase exactly by one at every step; this is the
case of maximal fixity $n-1$. 
\end{remark}

We now describe a characterization of $p$--groups of fixity $f<p$,
which appears in \cite{NOS}. 

\begin{theorem}\label{groups-fixityone}
Let $f$ be a non--negative integer and let $p$ be a prime number
greater than $f$. A non--abelian 
$p$--group $P$ has fixity $f$ if and only if
the following hold:

\begin{itemize}
\item $P$ has $p$--rank equal to $f+1$;

\item $P$ has cyclic center and
an abelian maximal subgroup\footnote{By an abelian maximal subgroup we
mean a maximal proper subgroup which is abelian.}.

\end{itemize}

Moreover, for any such group there exists a
faithful irreducible $P$--module of fixity $f$ and dimension $p$.
\end{theorem}

The $p$--groups of fixity $f<p$ have been explicitly described by
Conlon in \cite{Co}, and enumerated in \cite{NOS}. 
Using this one can verify for example
that for $p\ge 5$ and $n\ge 4$, there are exactly \emph{four}
non--abelian $p$--groups of order $p^n$ and fixity one. They can be listed as
follows in terms of generators and relations (see 
\cite{NOS}, page
228 and \cite{Huppert}, page 343):

\begin{itemize}
\item $<a,b~|~a^{p^{n-1}}=b^p=1, a^b=a^{1+p^{n-
2}}>$

\smallskip

\item $<a,x,y ~|~ a^{p^{n-3}}=[x,y],
a^{p^{n-2}}=[a,x]=[a,y]=x^p=y^p=1>$

\smallskip

\item $<a,x ~|~ a^{-\lambda p^{n-3}}=[x,a,x], a^{p^{n-2}}=x^p= [x,a]^p=[x,a,a]=1>$, where
$\lambda$ is equal to $1$ or to a non--quadratic residue modulo $p$.

\end{itemize}
\noindent Later we will construct 
actions of these $p$--groups on a product of two spheres. Note that the first
group on this list is a metacyclic group; the other three are said
to be of \emph{exceptional type} and they appear in the classification of rank
two $p$--groups which we will make use of in a subsequent section.
The situation for $p=2$ is rather different: there are $3n-8$ nonabelian
groups of order
$2^n$ having fixity equal to one.

In \cite{Sh} a detailed analysis of group theoretic properties
of groups of low fixity is carried out. 
The main result is that there exists a function
$\theta :\mathbb N\to\mathbb N$ such that if $V$ is a $\C G$--module of fixity
$f$, then there exists a normal subgroup $N\lhd G$ such that
$[G:N]\le \theta (f)$ and $N$ is solvable of derived length
at most $3$. Another result is that there exists a function
$\Phi :\mathbb N\to\mathbb N$ such that if $V$ is a $\C G$-module
of fixity $f$ and $\dim  V>\Phi(f)$, then either $G$ is
solvable or $G/O(G)$ has a subgroup of index at most 2 which is
isomorphic to $\rm{SL}(2,5)$. Here $O(G)$ denotes the largest
normal subgroup of $G$ whose order is odd, hence it is solvable.

There are substantial
restrictions for groups of low fixity. 
For example, if $G$ is a finite group acting on $V$ 
with fixity equal to one, 
then we have

\begin{itemize}

\item $\dim V\le |G|$. This is in contrast to the situation for
fixity zero, where we can take modules of arbitrarily large dimension. In fact
this is the only fixity where the modules do not have bounded dimension.

\item
The $p$--rank of $G$ is one for all primes $p< \dim V -1$. For example
this means that as soon as $\dim V>3$, the $2$--Sylow subgroup of $G$
must be either cyclic or generalized quaternion. This excludes all non--abelian
simple groups.

\item If $p=2,3$ do not divide $|G|$, then there exists at most
one prime $p$ with the property that $G$ has a non--abelian Sylow $p$--subgroup.

\end{itemize}

As we shall see later, groups of fixity one constitute a class
of non--periodic groups that act freely and smoothly
on a product of two spheres, and their
structure may shed some light on the general conjecture that any finite group
of rank equal to two acts freely on a product of two
spheres.  In this paper we prove this conjecture for $p$-groups with $p > 3$.

\section{Actions of Groups of Fixity One and Two} 
In this
section we begin by discussing the special case of groups having fixity equal to
one. As we have seen they will act freely on some $U(n)/U(n-2)$, which
happens to be a rather nice manifold
which we can approach using very direct
methods.

\begin{theorem} \label{fixityone}
A subgroup  $G\subset U(n)$  of
fixity  one acts freely and smoothly on $X=\mathbb S^{2n-1}\times \mathbb
S^{4n-5}$.
\end{theorem}

\begin{proof}
We have a fibration
$$
\mathbb S^{2n-
3}=U(n-1)/U(n-2)\to U(n)/U(n-2)\to U(n)/U(n-1) = \mathbb S^{2n-1},
$$
with
structure group $U(n-1)$. It can be identified with the sphere bundle of an
associated $(n-1)$-dimensional complex vector bundle $\xi$.    In the next
section we will see that the vector bundle $\xi$ is always non-trivial, but here
we will note that $\xi \oplus \xi$ is trivial since it is classified by an element
of $\pi_{2n-1}(BU(2n-2)) \cong \pi_{2n-2}(U(2n-2)) \cong\pi_{2n-2}(U)$,
which is zero by Bott periodicity. 

Note that $\xi$ is a $G$--vector bundle, and
$G$ acts freely away from the zero section.  By considering the pullback
diagram
$$
\begin{CD} E(\xi \oplus \xi) @>>> E(\xi)\times E(\xi) \\@VVV
@VVV \\\mathbb S^{2n-1}  @>\Delta >> \mathbb S^{2n-1} \times \mathbb
S^{2n-1},
\end{CD}
$$
we see that $\xi \oplus \xi$ is also a $G$--vector bundle
with action free away from the zero section.  We can take $X = S(\xi\oplus
\xi)$ (the associated sphere bundle) to complete the proof.
\end{proof}

\begin{corollary}If $G$ is any finite subgroup of $SU(3)$
then it will act freely and smoothly on $\mathbb S^5\times\mathbb
S^7$.
\end{corollary}

\begin{example}
The finite linear groups\footnote{To be precise we are referring
to primitive unimodular irreducible groups.} 
in low dimensions have been completely classified and listed 
by Blichtfeldt and others 
(see \cite{B}, \cite{Feit}).
For a fixed prime $p$ let $H_p$ denote the semidirect product of the
form $E_p\times_T\rm{SL}_2(\mathbb F_p)$, where 
$$
E_p =\langle a,b,c ~|~ a^p = b^p = c^p = 1, [a,b]=c\rangle
$$ is the 
extraspecial
$p$--group of order $p^3$ and exponent $p$. Note that $|H_p|=p^4(p-
1)(p+1)$.
We find four interesting subgroups of $SU(3)$: the alternating group
$A_5$, the simple linear group $SL_3(\mathbb F_2)\cong PSL_2(\mathbb
F_7)$, 
the triple cover $3A_6$, and $H_3$. We obtain

\begin{corollary}
The groups $A_5$, $SL_3(\mathbb F_2)$, $3A_6$, and
$H_3$ all act freely and smoothly on $\mathbb S^5\times\mathbb
S^7$. 
\end{corollary}
\end{example}

\begin{example}
In \cite{Sh} page 290, a group $G_p\subset U(p)$ of fixity equal
to one and having order $2p^3(p-1)$
is explicitly described as follows. 
Let $v_0,v_1,\dots , v_{p-1}$ denote a basis for the underlying
complex vector space and let $\omega$ denote a primitive
$p$-th root of unity. Take the two linear transformations
$A,B\in U(p)$ defined by $Av_i=\omega v_i$ and $Bv_i=\omega^iv_i$
where $0\le i\le p-1$. Now identify the set $\{0,1,\dots ,p-1\}$
with $\mathbb F_p$, and let $x$ be a generator of $\mathbb F_p^*$.
Let $\mu\in\mathbb C$ be a primitive $2(p-1)$-th root of unity.
We define $C,D\in U(p)$ via
$Cv_i=v_{i+1}$, $Dv_i=\mu v_{xi}$, where $i\in\mathbb F_p$.
Let $G_p$ denote the subgroup of $U(p)$ generated by $A,B,C,D$.
The subgroup $H=\langle A,B\rangle\cong\mathbb Z/p\times\mathbb Z/p$
is normal
in $G_p$, and $G_p/H$ is isomorphic to a double cover of the semi-direct product
$\mathbb Z/p\times_T\mathbb Z/(p-1)$. The group $G_p$ acts freely and
smoothly on $\mathbb S^{2p-1}\times\mathbb S^{4p-5}$. 
\end{example}

These techniques are particularly effective for constructing actions of
small $p$--groups.
\begin{proposition}\label{smallpgroups}
Let $P$ denote a $p$--group of rank equal to two and with
order $|P|\le p^4$. Then $P$ acts freely and smoothly on 
$\mathbb S^m\times\mathbb S^n$ for some $m,n>0$.
\end{proposition}
\begin{proof}
If the group $P$ has center $Z(P)$ of rank equal to two, then
every element of order $p$ in $P$ is central. In this case the
group acts freely on a product of two representations (see
\cite{AS}, page 422). Hence we can assume that $Z(P)$ is cyclic.
For any $p$--group of order $p^4$, there exists an abelian subgroup
of order $p^3$ (see \cite{Suz}, page 85), and so we can assume that $P$
has an abelian maximal subgroup. By Theorem~\ref{groups-fixityone}
we infer that $P$ has fixity equal to one, whence by Theorem~\ref{fixityone}
it will act freely and smoothly on a product of two spheres, 
completing the proof.
\end{proof}

As
we have seen, if $G\subset U(n)$ has fixity one, then it will act on $U(n)/U(n-
1)= \mathbb S^{2n-1}$ with rank one isotropy subgroups. In \cite{AS}, it was
shown that given a linear $G$--sphere $S(V)$ satisfying this condition, there
exists a finite $G$--CW complex $X$ with a free $G$--action, and such that
$X\simeq S(V)\times \mathbb S^M$ for some integer $M$
(see \cite{Unlu} for a more direct proof). From the above we
see that for a group $G$ of fixity one, this construction can be realized
explicitly via a free and smooth $G$--action on $\mathbb S^{4n-5}\times
\mathbb S^{2n-1}$.  However there are groups which act on a linear sphere
with periodic isotropy which do not have fixity equal to one. Indeed in
\cite{AS} it was shown that \emph{every} rank two $p$--group admits a
representation in $U(|P:Z(P)|)$ such that the action on $U(|P:Z(P)|)/U(|P:Z(P)|-
1)$ has periodic isotropy; on the other hand we have seen that $p$--groups
of fixity equal to one are rather restricted. However in the case of rank two
$p$--groups we have complete information:

\begin{proposition}\label{rank2}
If $p>3$ and $P$ is a finite $p$--group of rank equal to two, then
$P$ is either a metacyclic group or a group of exceptional type, with
fixity equal to one.
\end{proposition}
\noindent This proposition is proved in \cite{NOS}, page 228 by using the 
classification of rank two $p$--groups appearing in \cite{Blackburn}, as
well as the characterization of $p$--groups of fixity equal to
one described in Lemma~\ref{groups-fixityone}.
Next we recall the proof of an elementary
\begin{lemma}\label{metacyclic}
If $P$ is a metacyclic $p$--group then it acts freely on a product of
two linear spheres.
\end{lemma}
\begin{proof}
The group $P$ is an extension of the form $1\to \mathbb Z/p^t\to P
\to \mathbb Z/p^s\to 1$. Let $\chi$ denote a one dimensional character
of the subgroup $B\cong \mathbb Z/p^t$ which maps a generator of $B$
to a primitive $p^t$--th root of unity. Let $V = \ind^P_B(\chi)$; this is
a $p^s$--dimensional complex representation of $P$. The action of $P$
on the associated sphere $S(V)$ restricts to a free $B$--action; 
this can be checked using Mackey's formula (it suffices to show that the
unique and central subgroup of order $p$ in $B$ acts freely). Now let
the quotient group $A=P/B$ act freely on $\mathbb S^1$ through multiplication
by a primitive $p^s$-th root of unity. From this we can define a diagonal
$P$ action on $\mathbb S^{2p^s-1}\times\mathbb S^1$ which is evidently free.
\end{proof}

Obviously the action described above is a \emph{product} action. An interesting
fact is that for $p>3$, any exceptional rank two $p$--group does not admit a free
product action; indeed we have

\begin{proposition}
If $p>3$ and $P$ is a rank two $p$--group which is not metacyclic, then it
cannot act freely on $X=\mathbb S^m\times\mathbb S^n$ via a product action.
\end{proposition}
\begin{proof}
By the results in \cite{DH}, every action of a finite $p$--group on a sphere
can be modeled using a linear action; in particular given such a product action
there exist representations $V$, $W$ such that for any $p\in P$ we have
$X^{\langle p\rangle}\ne\emptyset$ if and only if $V^{\langle p\rangle}\ne\{0\}$
and $W^{\langle p\rangle}\ne\{0\}$. In particular this implies that 
if $P$ acts freely on $X$, then it acts freely
on $S(V)\times S(W)$. However it is shown in \cite{Ray}, page 486 that this is
impossible for a non--metacyclic $p$--group if $p>3$.
\end{proof}

We can now state and prove a geometric characterization of rank two
$p$--groups (for $p>3$) which naturally extends the classical rank one
situation, and which rather surprisingly has been a conjecture until
now:
\begin{theorem}
Let $p> 3$ be a prime number. Then a $p$--group $P$ acts freely 
and smoothly on $M=\mathbb S^m\times\mathbb S^n$ for some $m,n > 0$
if and only if
$P$ does not contain $\mathbb Z/p\times\mathbb Z/p\times\mathbb Z/p$
as a subgroup.
\end{theorem}
\begin{proof}
We have known for decades that $(\mathbb Z/p)^3$ does not act freely on a product of
two spheres (see \cite{Heller}). All we need to do is construct the actions.
By Proposition~\ref{rank2}, we know that $P$ is either metacyclic, or
of exceptional type and having
fixity equal to one. By Lemma~\ref{metacyclic}, the metacyclic case is
taken care of;
by Theorems~\ref{groups-fixityone} and \ref{fixityone},
the rank two $p$--groups of exceptional type will all act
freely and smoothly on $\mathbb S^{2p-1}\times\mathbb S^{4p-5}$.
\end{proof}

The situation at the `small primes' $p=2$ and $p=3$ is still unresolved.
It can be shown that for $n>4$, there exist exactly two $3$--groups of order
$3^n$ which are neither of fixity one, nor act freely on a product of two
linear spheres. At the prime $2$ the group theory is considerably
more complicated. We briefly summarize recent work by
\"O.\"Unl\"u in this direction (see \cite{Unlu}).
The smallest example of a rank two
$2$--group of fixity greater than one that does
not act freely on a product of two linear spheres is the extraspecial $2$--group 
$Q_8*D_8$, of 
order\footnote{By this we mean the central product 
$Q_8\times D_8 / \Delta$, where $\Delta\cong\mathbb Z/2$ is the central
diagonal subgroup.} $32$ (note that it does not contain a maximal subgroup
which is abelian). In fact one can show that a rank two $2$--group $P$ does
not contain $Q_8*D_8$ if and only if it either has fixity equal to one
or acts freely on a product of two linear spheres. Hence all such groups
act freely and smoothly on a product of two spheres. In order to handle this
particular group, we need the notion of \textsl{quaternionic fixity}, expressed
in terms of representations $G\subset Sp(n)$. In the case of $Q_8*D_8$, it 
embeds in $Sp(2)$, and twice the associated bundle is seen to split,
hence producing a free and smooth action of this group
on $\mathbb S^7\times\mathbb S^7$. Using these techniques, \"Unl\"u has shown
that among the $396$ groups of rank equal to two, and having order
which divides $256$, there is exactly \textbf{one} group which is not known to
act freely and smoothly on a product of two
spheres\footnote{Verifying this requires using the
computer algebra program MAGMA.}.

For fixity equal to two we can
apply the results in \cite{AS} to obtain a free action on a finite complex,
although not necessarily a manifold. As a consequence we will obtain new
examples of free actions for rank three groups. 

\begin{theorem}
If $G\subset
U(n)$ is of fixity equal to two, then $G$ acts freely on a finite complex
$X\simeq \mathbb S^{2n-1}\times\mathbb S^{4n-5}\times \mathbb S^M$
for some $M>0$. 
\end{theorem}

\begin{proof}
Consider the $G$--action on
$Y=U(n)/U(n-2)$; by our hypothesis it has rank one isotropy. Using the same
approach as in Theorem~\ref{fixityone} we obtain a $G$--action on
$Y'=\mathbb S^{2n-1}\times\mathbb S^{4n-5}$
with the same isotropy, hence of rank equal to one. Applying 1.4 in \cite{AS}
we obtain that $G$ acts freely on a finite complex $X\simeq
Y'\times\mathbb S^M$
for some $M>0$, which completes the proof.
\end{proof}

\begin{corollary}
If $G\subset SU(4)$ then $G$ acts freely on a finite complex
$X\simeq \mathbb S^7\times\mathbb S^{11}\times\mathbb S^M$ for some
$M>0$.
\end{corollary}

\begin{proof}
It suffices to observe that if $G\subset SU(4)$, then it will act on
$SU(4)/SU(2)$ with isotropy of rank at most one.
\end{proof}

\begin{example}
We are of course interested in groups of rank equal to three and of fixity
equal to two. The list of finite linear subgroups in \cite{Feit} yields two
interesting examples in $SU(4)$.  Let $T=\rm{SL}_2(\mathbb
F_3)*\rm{SL}_2(\mathbb F_3)=\rm{SL}_2(\mathbb
F_3)\times\rm{SL}_2(\mathbb F_3)/\Delta$, where $\Delta\cong\mathbb
Z/2$ is the diagonal subgroup of order two. Note that $|T| = 2^5\cdot 3^2$
and that its $2$--Sylow subgroup is the central product $Q_8*Q_8$, which is
an extra-special $2$--group of order $32$ and has rank equal to three. The second
example is $G=\rm{Sp}_4(\mathbb F_3)$, in this case
$|G|=2^6\cdot 3^4\cdot 5$. It has $2$--rank and $3$--rank both equal
to three.
Hence we obtain

\begin{corollary}
The groups $\rm{SL}_2(\mathbb F_3)*\rm{SL}_2(\mathbb F_3)$ and
$\rm{Sp}_4(\mathbb F_3)$ act freely on a finite complex
of the form 
$X\simeq \mathbb S^7\times\mathbb S^{11}\times\mathbb S^M$ for some
$M>0$.
\end{corollary}
\end{example}

\begin{example}
From the Atlas \cite{Con} we see that  
there is an embedding of $G=\rm{PSL}_2(\mathbb F_8)$ in $SU(7)$. 
Note that $|G|=2^3\cdot 3^2\cdot 7$
and that the $3$--Sylow subgroup is cyclic while the $2$--Sylow subgroup
$S$ is elementary abelian. From the character table in the Atlas we can
infer that this $7$--dimensional representation $V$ restricts to the reduced
regular representation on the subgroup $S$. Hence on any rank two
subgroup $E\subset S$, $V_{|_E}\cong I\oplus I\oplus \mathbb C$, where $I$ is
the $3$--dimensional reduced regular representation. This representation
clearly cannot be conjugated into the subgroup $SU(5)\subset SU(7)$ defined
as before by extending matrices in $SU(5)$ by the $2\times 2$ identity
matrix. Hence $E$ acts with cyclic isotropy on $SU(7)/SU(5)$. Applying the
techniques outlined above, we infer that $G=\rm{PSL}_2(\mathbb F_8)$ acts
freely on a finite complex $X\simeq \mathbb S^{13}\times\mathbb
S^{23}\times\mathbb S^M$for some $M>0$. Note that in this example $G$
acts freely on $SU(7)/SU(3)$ but not on $SU(7)/SU(4)$, hence $G\subset
SU(7)$ has fixity equal to three; nevertheless our methods can be applied to
this representation.
\end{example}

More generally if $G$ acts on $U(n)/U(n-f)$ with periodic isotropy,
the results in \cite{AS} establish the existence of a free action on  a finite
complex $X\simeq U(n)/U(n-f)\times\mathbb S^M$ for some large $M$.

In the next section we use surgery theory (propagation of
group actions) to construct a free action of a subgroup $G \subset U(n)$ of
fixity $f$ on a product of $f+1$ spheres provided the order of $G$ is prime
to $(n-1)!$.  If $G$ is such a group and has fixity one,  these techniques provide a
stronger result than Theorem \ref{fixityone} in the sense that the action is
on a lower-dimensional product of spheres. 

\section{Propagating Group Actions}
In this section we will show that  a free
action of a group of
order prime to $(n-1)!$ on a Stiefel manifold
$U(n)/U(k)$ propagates to a
free
action on a product of spheres.  To begin we recall the integral cohomology
of Stiefel manifolds.

\begin{proposition} \label{cohomology}
The Stiefel
manifold $U(n)/U(k)$ has the same integral cohomology as a product
$\mathbb S^{2n-1}\times\mathbb S^{2n-3}\times\dots\times \mathbb
S^{2k+1}$ of $n-k$ spheres.
\end{proposition}

\begin{proof}
We will prove this
by downward induction on $k$; clearly it is true for $k$ equal to $n-1$.  
Consider the fibration $\mathbb S^{2k+1}\to U(n)/U(k)\to U(n)/U(k+1)$; the
Euler class must be zero since, using the inductive hypothesis, it sits in a 
zero
group.  Hence the associated Gysin sequence breaks into a sequence of short
exact sequences.  We  infer that the cohomology of $U(n)/U(k)$ is an
exterior algebra on generators in the desired dimensions.
\end{proof}

One
would not expect the Stiefel manifold to have the same homotopy type as a
product of spheres, since Bott periodicity forces a regularity on the
homotopy groups of Stiefel manifolds which is absent from the homotopy
groups of a product of spheres. However, the differences only involve primes
less than $n$.  The next several results illustrate these
phenomena.

\begin{proposition} \label{Bott}
$\pi_{2n}(U(n))\cong \mathbb
Z/n!$ and the generator is given by the characteristic element of the
fibration
$U(n)\to U(n+1)\to S^{2n+1}$.
\end{proposition}

\begin{proof}
This result is
due to Bott \cite{Bo}, however we will sketch a proof which uses only stable
Bott periodicity. Note that $\pi_{2n}(U(n+1))=\pi_{2n}(U)=0$
and $\pi_{2n+1}(U(n+1))=\pi_{2n+1}(U)=\mathbb Z$ (these groups are in the
stable range), hence we have an exact
sequence
$$
\pi_{2n+1}(U(n+1))\xrightarrow{p_*}
\pi_{2n+1}(\mathbb S^{2n+1})\xrightarrow{\Delta_*} \pi_{2n}(U(n))\to 0.
$$
Recall (see
\cite{Wh}, page 206) that the characteristic element of the fibration above
is given by $\Delta_*(\iota_{2n+1})\in\pi_{2n}(U(n))$, where $\iota_{2n+1}$
is the canonical generator for $\pi_{2n+1}(\mathbb S^{2n+1})$.  Hence
$\pi_{2n}(U(n))$ is cyclic with generator the characteristic element. 

The
homomorphism $p_*$ can be identified with taking the Euler class of
the corresponding $(n+1)$-dimensional complex bundle over $\mathbb
S^{2n+2}$, and hence with the top Chern class $c_{n+1}$.  But the Chern
classes which arise from bundles over $\mathbb S^{2n+2}$ are precisely the
multiples of $n!$ (see \cite{Hu}, page 280). The result follows.
\end{proof}

The following technical proposition will be a key ingredient in
our construction
of group actions.

\begin{proposition} \label{product}
Given integers $n>k$,
there is a map
$$
\mathbb S^{2n-1}\times\mathbb S^{2n-3} \times\dots\times \mathbb S^{2k+1}\to
U(n)/U(k)
$$
which induces an isomorphism in homology with coefficients
in $\mathbb Z[\frac{1}{(n-1)!}]$.
\end{proposition}

\begin{proof}
Consider the
fiber bundle 
$$
U(n-1)/U(k)\to U(n)/U(k)\to\mathbb S^{2n-1}
$$
with
structure group $U(n-1)$.  By Proposition \ref{cohomology} (or a
Wang sequence) 
$$
H^*(U(n)/U(k);\mathbb Z) =  H^*(\mathbb S^{2n-1};\mathbb Z) \otimes
H^*(U(n-1)/U(k); \mathbb Z) .$$

By Proposition \ref{Bott} the fiber bundle is classified
by the homotopy class $\mathbb S^{2n-1}\to BU(n-1)$ representing the
generator of $\pi_{2n-1}(BU(n-1))= \mathbb Z/(n-1)!$. Hence if we take
a map $g:\mathbb S^{2n-1}\to\mathbb S^{2n-1}$ of degree $(n-1)!$ then the
induced bundle with fiber $U(n-1)/U(k)$ will be trivial, with total space
$E$
homeomorphic to $\mathbb S^{2n-1}\times U(n-1)/U(k)$; note that it comes
equipped with a map $E\to U(n)/U(k)$ which induces a $\mathbb
Z[\frac{1}{(n-1)!}]$ homology equivalence. Using downward induction on $n$
we can easily obtain the desired map.
\end{proof}

Hence the Stiefel manifold $U(n)/U(k)$ has the same
integral cohomology as a product of spheres, and after inverting
$(n-1)!~$, we can realize this isomorphism by a map from the product
of spheres to the Stiefel manifold. Under these conditions, if we
have a free action of a group on $U(n)/U(k)$ and if the order of the group is
prime to $(n-1)!~$, then the action will
``propagate'' to a free action on an actual product of spheres. 

\begin{theorem} \label{action}
Let $G$ denote a finite subgroup of $U(n)$ which
acts freely on
$U(n)/U(k)$ with $k\geq 1$. If
the order of $G$ is prime to $(n-1)!~$, then $G$ acts freely
and smoothly on $\mathbb S^{2n-1}\times\mathbb S^{2n-3} \times\dots\times \mathbb S^{2k+1}$.
\end{theorem}

\begin{remark}  In the proof of the theorem we use the fact that $U(n)/U(k)$ is
simply-connected when $k\geq 1$.  It
seems likely that Theorem \ref{action} is true in the non-simply-connected case 
$G \subset U(n)$ ($k = 0$), but the
surgery theoretic complications are considerable, and we will not consider them here.
However our methods will still allow us to construct a free $G$--action on a finite 
complex with the
homotopy type of a product of $n$ spheres 
(see Corollary \ref{non-simply-connected}). As we remarked earlier
(see \ref{SU(n)}), our theorem always applies to $G\subset SU(n)$, as
it acts freely on $U(n)/U(1)$.
\end{remark}

The proof of Theorem \ref{action}
uses the  method of
propagation of group actions.   We refer   to 
\cite{CW} for the 
basic  technique, to \cite{AD} for a short survey,  to
\cite{DW1} for a nice application, to
\cite{DL} for a key lemma, and to \cite{DW2} for definitive statements of
results.  The philosophy is that if two manifolds (e.g.
the Stiefel manifold and the product of spheres) resemble each other
homologically at the order of a
finite group $G$, then their behavior with
respect to $G$--actions should be similar.  To get actions on manifolds, the
methods are surgery theoretic
and the technical details can be
formidable.  

\begin{definition} 
A $G$--action on  $Y$ {\em propagates across a map $f : X \to Y$} if there
is a $G$--action on $X$ and  an equivariant map homotopic to $f$.  
\end{definition}

\begin{definition} 
A $G$--action on a space $X$ is {\em
homologically trivial} if the induced action on   $H_*(X; \Z[1/q])$ is
trivial,
where
$q$ is the order of $G$. 
\end{definition}

The main theorem involving
propagation is stated below.  We will spend the rest of the section defining
terms in its statement, outlining the proof, and applying it to the case of
interest (Theorem \ref{action}).  Most of the theorem below is due to
Cappell-Weinberger \cite{CW}, however  the  general case is due to  Davis-L\"offler
\cite{DL}. We
outline the proof since that will make application of the theorem easier,
and
also because the full statement is not easy to find in the
literature.

\begin{theorem}
\label{propagation}  Let $f : X \to Y$ be a map between simply-connected
spaces having the homotopy type of CW-complexes.  Let $G$ be a group of
order $q$ acting freely and homologically trivially on $Y$.  Consider
the following conditions:  

\begin{enumerate}

\item  $f$ is a $\Z_{(q)}$--equivalence:
$$
f_* : H_*(X; \Z_{(q)}) \xrightarrow{\simeq} H_*(Y;
\Z_{(q)}).
$$

\item  The Swan obstruction
vanishes:
$$
\sigma(\chi^{\text{tor}}(f)) = 0 \in \widetilde K_0(\Z G).
$$
\item
$X$ and $Y$ are closed smooth manifolds of dimension greater than four and
that the action of $G$ on $Y$ is smooth.

\item The normal invariant of the
$\Z_{(q)}$--local homotopy equivalence $f$ is in the image of the transfer
map:
$$
\nu_{(q)}(f) \in \im(p^* :[Y/G,F/O]_{(q)} \to
[Y,F/O]_{(q)}).
$$
\end{enumerate}

Then 

(a)  If (1) holds, there is a CW-complex $X'$ and a homotopy equivalence $h : X' \to X$ so that the 
$G$--action
on $Y$ propagates across $f \circ h$, with a cellular $G$--action on $X'$. 
Furthermore, the homotopy type of $X'/G$ is uniquely determined.

(b)  If (1)
holds and $Y/G$ and $X$ have the homotopy type of  finite-dimensional CW-
complexes, there is a finite-dimensional CW-complex $X'$ satisfying the
conclusion of (a).

(c)  If (1) and (2) hold and $Y/G$ and $X$ have the
homotopy type of  finite CW-complexes, there is a finite CW-complex $X'$
satisfying the conclusion of (a).

(d) If $q$ is odd and (1), (2), and (3) hold,
then $X'$ can be taken to be a closed smooth manifold.

(e) If $q$ is odd and
(1), (2), (3),  and (4) hold, then the $G$--action on $Y$ propagates across $f$,
with
a smooth action on $X$.
\end{theorem}


Before sketching a proof of Theorem~\ref{propagation} we recall some
basic background material.
The homotopy aspects of
propagation of group actions
depend on localizing topological spaces (see \cite{HMR},
\cite{BK1}, and
\cite{BK2} for details). Let $R \subset \Q$ be a subring of the rationals.  For
an abelian group $A$, let $A_R = A \otimes R$.  A homomorphism $A \to B$ is
an {\em $R$-equivalence} if the induced map $A_R \to B_R$ is
an isomorphism.  An abelian group $A$ is {\em $R$-local} if the map $A \to
A_R$ is an isomorphism.

We assume any space discussed below has the homotopy type of
a connected CW-complex.  A map $X \to Y$ is an {\em $R$-equivalence} if the
induced map on the fundamental group is an isomorphism and the induced
maps on higher homotopy groups $\pi_n, n > 1$ are  $R$-local equivalences.
Equivalently, the map is an isomorphism on the fundamental group and the
induced map on the homology of the universal covers $H_*(\widetilde X) \to
H_*(\widetilde Y)$ is an $R$-equivalence.  A space $X$ is $R$-local
if $\pi_n(X)$
is $R$-local for all $n > 1$, equivalently $H_*(\widetilde X)$ is $R$-local.
An
{\em $R$-localization of $X$} is a $R$-equivalence $X\to Y$ where $Y$ is
$R$-local.  There are existence and uniqueness theorems for $R$-localizations of
$X$.  Their existence follows by applying  the fiberwise localization
theorem of
Bousfield-Kan \cite[p. 40]{BK2}  to $X \to B\pi_1X$.  Their $R$-localization
is
functorial on the geometric realizations of simplicial sets and maps.  If
$f : X\to Y$ and $g : X \to Z$ are two $R$-localizations of $X$, then there
is a homotopy equivalence $h : Y \to Z$ so that $h\circ f \simeq g$. We will
write
$X \to X_R$ to denote an $R$-localization of $X$. 

Let $q$ be a nonzero
integer.  Let $\Z_{(q)} = \Z[\frac{1}{p_1},\frac{1}{p_2}, \frac{1}{p_3},
\dots]\subset \Q$, where $\{p_1,p_2,p_3,\dots\}$ is the set of primes
which do {\em not} divide $q$. For a space $X$, let $X \to X_{(q)}$ and $X\to
X{[1/q]}$ denote the $R$-localizations of $X$ where $R$ equals $\Z_{(q)}$ and
$\Z[1/q]$ respectively.  Then $X$ is the homotopy pullback  of
$$
\begin{CD}@.
X{[1/q]}\\@. @VVV\\X_{(q)} @>>> X_\Q,\end{CD}
$$
(i.e. $X$ is homotopy
equivalent to  what results after converting the vertical map to a
fibration and taking the pullback). 

S. Weinberger \cite{We} made the
following
key observation.

\begin{lemma} A free $G$--action on a simply-connected 
CW-complex $X$ is homologically trivial if and only if $(X/G)[1/q] \simeq
X[1/q] 
\times BG$. 
\end{lemma}

\begin{proof}  Suppose a finite group $G$ of order
$q$ acts freely on a simply-connected space $X$.  Then $X \subset Y = X\cup
e^2_1 \cup \dots \cup e^2_k$ with $Y$ simply-connected (add 2-cells to kill
the fundamental group).  Since $H_1(X/G)= G/[G,G]$ is $q$-torsion,
$\pi_2(Y) \otimes \Z[1/q] \cong H_2(Y; \Z[1/q])\cong \Z[1/q]^k  $.  Let $Z
= Y \cup e^3_1 \cup\dots \cup e^3_k$  where the attaching maps of the 
3-cells represent a $\Z[1/q]$-basis of $\pi_2(Y) \otimes \Z[1/q]$.  By means
of this ``plus" construction we have constructed a map $i: X/G \to Z$ to
a simply-connected space inducing an isomorphism
on $H_*(\quad;\Z[1/q])$.

Now suppose, in addition, that $G$ acts
homologically trivially on $X$.  Then a transfer argument shows that the
covering map $\pi: X \to X/G$ induces an isomorphism on $H_*(~;\Z[1/q])$,
hence so does $i \circ \pi : X \to Z$.  Thus $Z[1/q]\simeq X[1/q]$, so  we
have  $\Z[1/q]$-equivalence $X/G \to  X[1/q]  \times BG$ to a 
$\Z[1/q]$-space.  By uniqueness of localization, there is a
homotopy equivalence
$(X/G)[1/q] \to X[1/q] \times BG $ as desired.  

Conversely if $(X/G)[1/q]
\simeq X[1/q] \times  BG $, then the $G$--action on $X$ is homologically
trivial
since the $G$--action on $X[1/q]  \times EG  $ clearly is.
\end{proof}

One
could
also prove the above lemma by using obstruction theory to show that the
fibration
$$
X[1/q] \to X[1/q]\times_G EG \to BG 
$$
is fiber homotopically
trivial and noting $  X[1/q]\times_G EG \simeq X/G[1/q]$.

\begin{proof}[Proof of Theorem \ref{propagation} (a)]  Let $X'/G$ be a
CW-complex having the
homotopy type  of the homotopy pullback of 
$$
\begin{CD}
@. X{[1/q]}\times
BG\\@. @VVV\\(Y/G)_{(q)} @>>> Y_\Q \times BG
\end{CD}
$$
where the
vertical map is given by applying $f[1/q]$ and then  $\Q$-localization, and the
horizontal map is provided by applying  $\Q$-localization and then
Weinberger's Lemma.
\end{proof}

\begin{proof}[Proof of Theorem \ref{propagation} (b)]  
We will
use the criterion of Wall
\cite[Thm. E]{Wa1} which says that a CW-complex $Z$ has the homotopy type
of an complex of dimension $N$ if and only if $H^i(\widetilde Z; \mathbb Z) = 0$ for 
all $i
> N$ and $H^{N+1}(Z;M)=0$ for all local coefficient systems $M$.  

Recall
that
both $Y/G$ and $X$ are assumed to have the homotopy type of 
finite dimensional CW-complexes; let $N$ be greater than or equal to
 both dimensions. Let $X'/G$ be a CW-complex produced by
Theorem \ref{propagation} (a). Then for $i > N$, $H^i(\widetilde{X'/G}; \mathbb Z) =
H^i(X; \mathbb Z) = 0$. Let $M$ be a local coefficient system for $X'/G$. For a prime 
$p$
dividing the order of $G$, 
$H^{N+1}(X'/G;M)_{(p)}$ is isomorphic to
$$H^{N+1}((X'/G)_{(p)};M_{(p)}) \cong H^{N+1}((Y/G)_{(p)};M_{(p)})
\cong H^{N+1}(Y/G; M)_{(p)}=0.$$
Likewise for a prime $p$ not
dividing the order of $G$, a transfer argument
shows
$$H^{N+1}(X'/G;M)_{(p)} \cong H^{N+1}((X'/G)_{(p)};M_{(p)}) 
\cong H^{N+1}(X'_{(p)};M_{(p)})
\cong H^{N+1}(X; M)_{(p)}=0.$$
Hence $H^{N+1}(X'/G;M) = 0$.  Thus $X'/G$ has the homotopy type of an
$N$-complex.
\end{proof}

This leads to a nontrivial result.

\begin{corollary}  Let $G$ denote a finite subgroup of $U(n)$ which
acts freely on
$U(n)/U(k)$ with $k\geq 1$.  Suppose the  order of $G$ is prime to $(n-1)!~$. 
Then $G$ acts freely on a finite dimensional CW-complex having the homotopy
type of  $\mathbb S^{2n-1}\times \mathbb S^{2n-3}
\times\dots\times \mathbb S^{2k+1}$.
\end{corollary}

It is traditional in the study of actions on products of spheres to be satisfied
with an action on a finite-dimensional
complex which has the homotopy type of a product of spheres.  But that is
not good enough for us; we won't stop until we
have constructed an action on the product of spheres themselves!  Next we
will try and improve from a finite-dimensional
complex to a finite complex.


Wall \cite{Wa1} showed
that a connected CW-complex $W$ has the homotopy type of a finite 
CW-complex if and only if $\pi_1(W)$ is finitely presented and
$C_*(\widetilde W)$ is
chain homotopy equivalent to finite free $\Z[\pi_1W]$-chain complex.
If $C_*(\widetilde W)$ is finitely dominated this is equivalent to the vanishing
of the finiteness obstruction $[C_*(\widetilde W)]\in \widetilde K_0(\mathbb Z[\pi_1(W)])$.

We recall the definition of the Swan homomorphism appearing
in Theorem~\ref{propagation}.
For an integer $r$ relatively prime to $q=|G|$,
define the projective $\Z G$-module $P_r =\ker \epsilon : \Z G \to \Z/r$. 
One can show that
$P_r\oplus P_s
\cong P_{rs}
\oplus \Z G$ and that $P_{1+kq}\cong \Z G$.
Hence there is a well-defined \emph{Swan homomorphism}
$\sigma : (\Z/q)^\times \to \widetilde
K_0(\Z G)$, defined by $\sigma([r]) = [P_r]$.
Now given a rational homology equivalence $f : X \to Y$, define $\chi^{\tor}(f)$ to
be the rational number  
$
\prod_i |H_i(C_f)|^{(-1)^i}$, where $C_f$ is the mapping cone of $f$.
If $f$ is a $\Z_{(q)}$-equivalence, then $\chi^{\tor}(f) \in
\Z_{(q)}^\times$ and hence determines an element
$[\chi^{\tor}(f)] \in (\Z/q)^\times$. 

\begin{proof}[Proof of Theorem \ref{propagation}(c)]  
By Theorem
\ref{propagation} (a) there is a CW-complex $Z$, a homotopy equivalence
$h_Z : Z
\to X$, a free, cellular $G$--action on $Z$, and an equivariant map $f_Z : Z \to
Y$ so that $f_Z \simeq f \circ h_Z$.  There is
a short exact sequence of chain complexes
$$
0 \to C_*(Y) \to C_*(C_{f_Z}) \to \Sigma C_*(Z) \to 0.$$
Our goal is to show that $C_*(Z)$ is finitely dominated and that the finiteness
obstruction $[C_*(Z)]\in\widetilde K_0(\mathbb ZG)$ vanishes.
Since $\oplus H_i(C_{f_Z})$ is finitely generated and each $H_i(C_{f_Z})$ has a finite
resolution over $\Z G$, $C_*(C_{f_Z})$ is
finitely dominated, as is $C_*(Y)$.  Hence $\Sigma C_*(Z)$ is finitely
dominated and
$$[C_*(Z)] = - [\Sigma C_*(Z)]
= [C_*(Y)] - [C_*(C_{f_Z})]
=- \sigma ([\chi^{\tor}(f)]).$$
Thus $C_*(Z)$ is chain homotopy equivalent to a finite free $\Z G$ chain
complex (equivalently $Z/G$ has the homotopy type
of a finite CW-complex) if and only if $\sigma([\chi^{\tor}(f)]) = 0$.
\end{proof}

\begin{corollary}  
Let $G$ be a subgroup of $U(n)$ which acts freely on
$U(n)/U(k)$  with $k\geq 1$ and whose order is prime to $(n-1)!~$. 
Then $G$ acts freely on a finite  CW-complex having the homotopy type of 
$\mathbb S^{2n-1}\times \mathbb S^{2n-3}
\times\dots\times \mathbb S^{2k+1}$.
\end{corollary}

\begin{proof}  
Here are two properties of $\chi^{\tor}$, whose verification
is left to the reader.

\begin{itemize}
\item $\chi^{\tor}(f \circ g) = \chi^{\tor}(f)\chi^{\tor}(g)$.
\item If $A_*$, $B_*$, and $C_*$ are graded abelian groups with $C_*$ free,
and if $f_* : A_* \to B_*$ is a graded map, then $
\chi^{\tor}(f_* \otimes \id_{C_*}) = \chi^{\tor}(f_*)^{\chi(C_*)}.
$
\end{itemize}

In the proof of Proposition \ref{product} we constructed a 
$\Z_{(q)}$-equivalence 
$$
f_{n,k} : \mathbb S^{2n-1} \times U(n-1)/U(k) \to U(n)/U(k)
$$
so that on homology
$${f_{n,k}}_*= g_* \otimes \id_* : H_*(\mathbb S^{2n-1}) \otimes H_*(U(n-
1)/U(k)) \to H_*(\mathbb S^{2n-1}) \otimes H_*(U(n-1)/U(k)), $$
where $g : \mathbb S^{2n-1} \to \mathbb S^{2n-1}$ is a map of degree $(n-
1)!~$.  Thus $\chi^{\tor}({f_{n,k}}) = ((n-1)!)^0 = 1.$

The $\Z_{(q)}$-equivalence produced by Proposition \ref{product}
$$
f : \mathbb S^{2n-1}\times\dots\times \mathbb S^{2k+1}\to U(n)/U(k)
$$
is the composite
$$
f = (\id_{\mathbb S^{2n-1} \times \dots \times \mathbb S^{2k+5}} \times
f_{k+2,k}) \circ \dots \circ (\id_{\mathbb S^{2n-1}} \times f_{n-1,k})
\circ f_{n,k}
$$
and hence $\chi^{\tor}(f)
= 1$, and so the result follows from  Theorem \ref{propagation}(c).
\end{proof}

\begin{corollary}  Let $G$ be a subgroup of $U(n)$ which acts  on
$U(n)/U(k)$  with $k>0$, having rank one isotropy and whose order is prime 
to $(n-1)!~$. 
Then $G$ acts freely on a finite  CW-complex having the homotopy type of 
$\mathbb S^{2n-1}\times \mathbb S^{2n-3}
\times\dots\times \mathbb S^{2k+1}\times \mathbb S^{M}$ for some $M > 1$.
\end{corollary}

\begin{proof}  By \cite{AS}, $G$ acts freely on a finite complex $X$ having the 
homotopy type of $U(n)/U(k) \times \mathbb S^{M}$ for some
$M > 1$.  One then propagates across the map 
$$
f \times \id : \mathbb S^{2n-1}\times \mathbb S^{2n-3}
\times\dots\times \mathbb S^{2k+1} \times \mathbb S^{M}\to U(n)/U(k) \times \mathbb S^{M} \simeq X.
$$
\end{proof}

\noindent The main application is the following intermediate result in the
non--simply connected case:
\begin{corollary}\label{non-simply-connected}
If $G\subset U(n)$ and $|G|$ is prime to $(n-1)!$, then $G$ acts freely on a
finite complex $X\simeq \mathbb S^{2n-1}\times\dots\times\mathbb S^3\times 
\mathbb S^M$ for some $M > 1$.
\end{corollary}
\begin{proof}
It suffices to observe that $G$ will act on $U(n)/U(1)$ with periodic isotropy,
whence we can apply the previous corollary.
\end{proof}

We do not give the outline of the surgery theoretic proof
of Theorem \ref{propagation}(d) and (e), but instead refer
to
the original source \cite{CW} which proves the theorem for $p$--groups and to \cite{DL} 
which proves the key fact needed to prove the theorem for general groups $G$.

We are now ready to prove our main theorem by checking the
propagation hypotheses. 

\begin{proof}[Proof of Theorem \ref{action}]
Let $X = \mathbb S^{2n-1}\times \mathbb S^{2n-3}
\times\dots\times \mathbb S^{2k+1}$ and $Y = U(n)/U(k)$ and let $f : X \to
Y$ be the map produced by Proposition
\ref{product}.  We have seen that conditions (1), (2), and (3) of Theorem
\ref{action} hold; we only need to verify condition
(4).  

The definition of the local
normal invariant $\nu_{(q)}(f) \in [Y, F/O]_{(q)}$ is quite subtle (see 
\cite[p. 12]{DW2} for a definition), but we will sidestep the subtleties by
noting that
 $[Y,F/O]_{(q)} \to [Y,F/O]_{(0)}$ is injective, since $H_*(Y; \mathbb Z)$ is torsion-free.
Since
$F/O  \to BO$ is a rational equivalence \cite{MM},
it suffices to show that the 
image of the local normal invariant in $[Y,BO]_{(0)}$ is trivial. Its image in
$[Y,BO]_{(0)}$ is given by the
difference of $((\deg f) (f_{(0)}^{-1})^*\tau_X)
-
\tau_Y$.  But the tangent bundle $\tau_X$ is stably trivial, and, according to
\cite{Su},
$\tau_Y$ is trivial.
\end{proof}

Given a group of rank equal to $r$, a challenging open problem
is to construct a free action of $G$ on a finite complex with the
homotopy type of a product of $r$ spheres (see \cite{AS}). Constructing
free actions on actual products of spheres is of course the final goal,
and progress on this beyond spherical space forms has been very
scant.\footnote{Indeed, we still cannot verify if every 
finite group acts freely and
homologically trivially on some product of spheres.}
The following application of Theorem~\ref{action} 
provides an infinite number of new examples.

\begin{corollary}
Let $P$ denote a finite non--abelian $p$--group with cyclic center
and an abelian maximal subgroup. If the rank of $P$ is $r$ and
$r<p$, then there exists a free, smooth and homologically trivial action
of $P$ on $M=\mathbb S^{2p-1}\times\dots\times\mathbb S^{2(p-r)+1}$,
a product of $r$ spheres.
\end{corollary} 

\begin{proof}
Given a group of this type, we know from the results in \S 2 that
there is a faithful irreducible representation $P\subset U(p)$
of fixity equal to $r-1$. Hence $P$ acts freely on $U(p)/U(p-r)$. Using the
propagation theorem, we infer that in fact $P$ must act freely and smoothly
on the
stated product of spheres.
\end{proof}

\begin{example}
Let $P$ denote any $p$--group of fixity equal to one, where $p$ is an
odd prime. 
Then $P$ has a faithful
$p$--dimensional
representation such that $P$ acts freely on $U(p)/U(p-2)$;
hence it will act freely and smoothly
on $M=\mathbb S^{2p-1}\times\mathbb S^{2p-
3}$. Note that this provides an action on a lower dimensional manifold
than that provided by Theorem~\ref{fixityone}. 
\end{example}

\begin{example}
Let $P=\mathbb Z/p\wr\mathbb Z/p$ (the wreath product); 
this is a group of order $p^{p+1}$, with
cyclic center and having a maximal subgroup which is abelian, $(\mathbb Z/p)^p$.
Indeed, there is a natural embedding $P\subset U(p)$, where the elementary
abelian subgroup is mapped to the diagonal matrices and the $\mathbb Z/p$ action
is represented by a permutation matrix. Applying Corollary~\ref{non-simply-connected},
we obtain that
$P=\mathbb Z/p\wr\mathbb Z/p$ acts freely on a finite complex with the
homotopy type of a product of $p$ spheres, namely
$X\simeq\mathbb S^{2p-1}\times\dots\times\mathbb S^3\times\mathbb S^M$ for
some $M > 1$.
\end{example}

\end{document}